\documentclass[12pt,letterpaper,dvips]{article}

\usepackage{amsmath,amsfonts,amssymb,amsthm}
\usepackage[all]{xy}

\swapnumbers
\theoremstyle{plain}
\newtheorem{lemma}{Lemma}[section]
\newtheorem{prop}[lemma]{Proposition}
\newtheorem{thm}[lemma]{Theorem}

\newtheorem{cor}[lemma]{Corollary}

\newtheorem*{claim*}{Claim}

\newtheorem{curves_excl_thm}[lemma]{Curve Exclusion Theorem}

\theoremstyle{definition}
\newtheorem{defn}[lemma]{Definition}
\newtheorem{defns}[lemma]{Definitions}
\newtheorem{ex}[lemma]{Example}

\newtheorem{emp}[lemma]{}

\newtheorem{prob}[lemma]{Problem}
\newtheorem*{rk*}{Remark}
\newtheorem*{rks*}{Remarks}

\newtheorem{famous95}[lemma]{The famous 95 families}
\newtheorem{proofofcthm31s}[lemma]{Proof of Theorem~\ref{thm:curves}
assuming $\mathbf{a_0 = a_1 = a_2 = 1}$}
\newtheorem{proofofcthm11}[lemma]{Proof of Theorem~\ref{thm:curves}
assuming $\mathbf{a_1 > 1}$}
\newtheorem{proofofcthm21s}[lemma]{Proof of Theorem~\ref{thm:curves}
assuming $\mathbf{a_1 = 1}$, $\mathbf{a_2 > 1}$}

\theoremstyle{remark}
\newtheorem*{notn*}{Notation}

\newcommand{\HH}{\mathcal H}
\newcommand{\II}{\mathcal I}
\newcommand{\LL}{\mathcal L}

\newcommand{\NNN}{\mathcal N}
\newcommand{\PP}{\mathbb P} 
\newcommand{\PPP}{\mathcal P}
\newcommand{\QQ}{\mathbb Q}

\newcommand{\ZZ}{\mathbb Z} 

\renewcommand{\le}{\leqslant}
\renewcommand{\ge}{\geqslant}

\newcommand{\dashto}{\dashrightarrow}

\newcommand{\qeq}{\sim_\QQ}
\newcommand{\iso}{\simeq}

\newcommand{\newmid}{\!\mid\!}
\newcommand{\Oh}{\mathcal O}

\newcommand{\Bs}{\operatorname{Bs}}
\newcommand{\Centre}{\operatorname{Centre}}
\newcommand{\Cl}{\operatorname{Cl}}
\newcommand{\CS}{\operatorname{CS}}
\newcommand{\CSXnH}{\operatorname{CS}\left(X,\textstyle\frac{1}{n}\HH\right)}
\newcommand{\Diff}{\operatorname{Diff}}

\newcommand{\hcf}{\operatorname{hcf}}

\newcommand{\KXnH}{K_X + \textstyle\frac{1}{n}\HH}

\newcommand{\mult}{\operatorname{mult}}
\newcommand{\None}{\operatorname{N^1}}
\newcommand{\NEbar}{\operatorname{\overline{NE}}}
\newcommand{\Nonsing}{\operatorname{Nonsing}}

\newcommand{\VXnH}{\operatorname{V_0}\left(X,\textstyle\frac{1}{n}\HH\right)}
\newcommand{\Vzero}{\operatorname{V_0}}
\newcommand{\XnH}{\left(X,\textstyle\frac{1}{n}\HH\right)}

\begin{document}

\title{The Curve Exclusion Theorem for elliptic \\ and K3~fibrations
  birational to \\ Fano 3-fold hypersurfaces}
\author{Daniel Ryder}
\date{June~2006}
\maketitle

\begin{abstract}
The theorem referred to in the title is a technical result that is needed
for the classification of elliptic and K3~fibrations birational to Fano
3-fold hypersurfaces in weighted projective space.  We present a complete
proof of the Curve Exclusion Theorem, which first appeared in the author's
unpublished PhD~thesis~\cite{Ry02} and has since been relied upon in
solutions to many cases of the fibration classification problem.
We give examples of these solutions and discuss them briefly.
\end{abstract}

\section{Introduction} \label{sec:intro}
The problem that motivates the work presented here is the following.

\begin{prob} \label{prob:main}
Let~$\:X = X_d \subset \PP(1,a_1,a_2,a_3,a_4)$ be a Fano 3-fold weighted
hypersurface in one of the `famous 95' families of Fletcher and
Reid~\cite{Fl00}.  Assuming that $X$~is general in its family, we seek to
classify the set of K3~fibrations $g \colon Z \to T$ with $Z$~birational
to~$X$ and the set of elliptic fibrations $g \colon Z \to T$ with
$Z$~birational to~$X$.
\end{prob}

Solutions to both the K3 and elliptic cases of this problem for families~1
and~3 of the~95 first appeared in papers of Cheltsov (see~\cite{Ch00},
\cite{Ch03} and further references therein).  These are the only two
of the 95~families whose members are smooth: $X = X_4 \subset \PP^4$ in
family~1 is a smooth quartic 3-fold and $X = X_6 \subset \PP(1,1,1,1,3)$ in
family~3 is a double cover of~$\PP^3$ branched in a smooth sextic.  For
four of the 93~remaining singular families, solutions to both the~K3 and
elliptic cases of Problem~\ref{prob:main} appeared in~\cite{Ry06} and one
other case, family~5, was dealt with earlier in the
unpublished~\cite{Ry02}.  Here is an example solution.

\begin{thm}[\cite{Ry06}]
 \label{thm:75main}
Let $X = X_{30} \subset \PP(1,4,5,6,15)_{x,\,y,\,z,\,t,\,u}$ be a general
member of family~75 of the~95.
\begin{itemize}
\item[(a)] Suppose $\Phi \colon X \dashto Z/T$ is a birational map from $X$
  to a K3 fibration $g \colon Z \to T$ (see~\ref{defns:fibr} below for our
  assumptions on K3 fibrations, and also on elliptic fibrations). Then
  there exists an isomorphism $\PP^1 \to T$ such that the diagram below
  commutes, where $\pi = (x^4,y) \colon X \dashto \PP^1$.
\[ \xymatrix@C=1.6cm{
X \ar@{-->}[r]^{\Phi} \ar@{-->}[d]_{\pi} & Z \ar[d]^g \\
\PP^1 \ar[r]^{\iso} & T \\
} \]
\item[(b)] There does not exist an elliptic fibration birational to $X$.
\item[(c)] If $\Phi \colon X \dashrightarrow Z$ is a birational map from
$X$ to a Fano 3-fold $Z$ with canonical singularities then $\Phi$ is
actually an isomorphism (so in particular $Z \iso X$ has terminal
singularities).
\end{itemize}
\end{thm}

The proof of this theorem relies on one particular case of our Curves
Exclusion Theorem~(\ref{thm:curves}~below); \cite{Ry06}~contains a proof of
this case, but no others.

In~\cite{Ch05} Cheltsov, building on previous joint work with
Park~\cite{CP} and on~\cite{Ry02}, was able to classify elliptic fibrations
birational to a general member of any of the 95~families, i.e., to solve
completely the elliptic case of Problem~\ref{prob:main}.  Both~\cite{CP}
and~\cite{Ch05} rely on Theorem~\ref{thm:curves}: see below.  One important
observation in these two papers --- which also appears in a simple form
in~\cite{Ch00} --- is that surprisingly useful information can be extracted
from the trivial fact that, in the elliptic case, the linear system on~$X$
with which we are working is not a pencil (see, e.g., Lemma~2.1
of~\cite{CP} and the proof of Lemma~2.11); largely because of this
observation, these papers deal only with the elliptic case of the
classification problem.  It should be noted, though, that~\cite{CP},
building on~\cite{Ry02}, contains constructions of K3~fibrations birational
to general members of all 95~families: it is the problem of excluding other
possible K3 fibrations that remains open, for the moment, in most cases.

Here is a theorem from~\cite{CP} which relies on
our~Theorem~\ref{thm:curves}.

\begin{thm}[{\cite[1.2]{CP}}] \label{thm:ell_fibrs}
A general $X_d \subset \PP(1,a_1,a_2,a_3,a_4)$ in family~$N$ of the~95 is
birational to an elliptic fibration if and only if \[ N \not\in
\{3,\,60,\,75,\,84,\,87,\,93\}. \]
\end{thm}

Theorem~\ref{thm:curves} is used in the proof of this result to help
demonstrate the nonexistence of a birational elliptic fibration for $N \in
\{3,\,60,\,75,\,84,\,87,\,93\}$.  Similarly, our theorem is used
throughout~\cite{Ch05} (see his Theorem~1.15 and Lemma~1.16) to classify
elliptic fibrations birational to all the 95~families.  We give one
example:

\begin{thm}[{\cite[26.3]{Ch05}}] \label{thm:36}
Let $X = X_{18} \subset \PP(1,1,4,6,7)_{x_0,\,x_1,\,y,\,z,\,t}$ be a
general member of family~36 of the~95 and assume that $\Phi \colon X
\dashto Z$ is a birational map from~$X$ to an elliptic fibration $g \colon
Z \to T$.  Then either there exists a birational map $\PP(1,1,4) \dashto T$
such that the diagram
\[ \xymatrix@C=1.6cm{
X \ar@{-->}[r]^{\Phi} \ar@{-->}[d]_{\pi} & Z \ar[d]^g \\
\PP(1,1,4) \ar@{-->}[r] & T \\
} \]
commutes, where $\pi = (x_0,x_1,y) \colon X \dashto \PP(1,1,4)$ is the
natural projection, or there exists a birational map $\PP(1,1,6) \dashto T$
such that the diagram
\[ \xymatrix@C=1.6cm{
X \ar@{-->}[r]^{\Phi} \ar@{-->}[d]_{\pi'} & Z \ar[d]^g \\
\PP(1,1,6) \ar@{-->}[r] & T \\
} \]
commutes, where $\pi' = (x_0,x_1,z) \colon X \dashto \PP(1,1,6)$ is the
natural projection.
\end{thm}

It is time to state the Curve Exclusion Theorem; first we need the
following.

\begin{notn*}
Let $X$ be a normal complex projective variety, $\HH$~a mobile linear
system on~$X$ and~$\,\alpha \in \QQ_{\ge\, 0}$.  We denote
by~$\CS(X,\alpha\HH)$ the set of centres on~$X$ of valuations that are
strictly canonical or worse for $K_X + \alpha\HH$ --- that is,
\[ \CS(X,\alpha\HH) = \{\Centre_X(E) \mid a(E,X,\alpha\HH) \le 0\}. \]
This notation is standard.  We also use the following nonstandard notation:
if $K_X + \alpha\HH$ is canonical then $\Vzero(X, \alpha\HH)$~denotes the
set of valuations (or of the corresponding divisors, each on some
sufficiently blown up model) which  are strictly canonical for~$K_X +
\alpha\HH$.
\end{notn*}

\begin{curves_excl_thm} \label{thm:curves}
Let~$\:X \,=\, X_d \,\subset\, \PP(1,a_1,a_2,a_3,a_4)_{x_0,\,\ldots,\,x_4}$
be a general hypersurface in one of the 95~families and $C \subset X$ a
reduced, irreducible curve. Suppose $\HH$~is a mobile linear system of
degree~$n$ on~$X$ such that $K_X + \frac{1}{n}\HH$ is strictly canonical
and $C \in \CSXnH$. Then there exists a pair~$\,\ell,\,\ell'\,$~of linearly
independent forms of degree~$1$ in~$(x_0,\ldots,x_4)$ such
that
\begin{equation} C \;\subset\; \{\ell = \ell' = 0\} \,\cap\,
  X. \label{CinPi}
\end{equation}
\end{curves_excl_thm}

\begin{emp}
For a precise discussion of how this theorem is used in the proofs of
Theorems~\ref{thm:75main}, \ref{thm:ell_fibrs}~and~\ref{thm:36} we
refer to the papers already cited; but we give a brief outline here.
Suppose we have a birational map $\Phi \colon X \dashto Z/T$ from a Fano
3-fold hypersurface $X = X_d \subset \PP(1,a_1,a_2,a_3,a_4)$ in one of the
95~families to either an elliptic or a K3~fibration $g \colon Z \to T$.  By
an analogue of the Noether--Fano--Iskovskikh inequalities --- which are
used in the Sarkisov program to break up a birational map between two Mori
fibre spaces into elementary links (see~\cite{Co95}) --- the log
pair~$\XnH$ has nonterminal singularities, where $\HH = \Phi^{-1}_* g^*
|A_T|$ is the transform on~$X$ of a very ample complete linear
system~$|A_T|$ on~$T$ and $n = \deg \HH$ is its anticanonical degree, i.e.,
$\HH \subset |{-n K_X}|$.  Using the main theorem of~\cite{CPR}, which
states that $X = X_d$~is \emph{birationally rigid}, we reduce to the case
where $\XnH$~has canonical but nonterminal, i.e., strictly canonical,
singularities.

At this point it is natural to ask what~$\CSXnH$ is; one of
the main results we use to answer this is our Curve Exclusion
Theorem~\ref{thm:curves}, which tells us that the only curves that could be
in~$\CSXnH$ are the obvious ones (see below).  Of course we also need
results describing which nonsingular and singular points could belong
to~$\CSXnH$, but we do not discuss this here.  Finally, given a complete
list of possibilities for~$\CSXnH$, we use various techniques to try to
deduce a complete list of birational elliptic and K3~fibrations (this is
a simplification of the process, but it gives the general~idea).

We expand a little on why it is `obvious' that certain curves cannot be
excluded.  We need the following:
\end{emp}

\begin{prop}[{\cite[2.2]{Ry02}}]
Let $X_d \subset \PP(1,1,a_2,a_3,a_4)$ be general in one of the families
with~$\,a_1 = 1$ and $\ell,\ell' \in k[x_0,\ldots,x_4]$ two independent
forms of degree~$1$. Then a general fibre $S$ of $\pi = (\ell,\ell') \colon
X \dashrightarrow \PP^1$ is a quasismooth Du Val K3 surface and, setting
$\PPP = \pi^{-1}_*\left|\Oh_{\PP^1}(1)\right|$,
\[ \CS(X,\PPP) = \{C_0,\ldots,C_r,P_1,\ldots,P_s\}, \]
where $C_0,\ldots,C_r$ are the components of $\{\ell = \ell' = 0\} \,\cap\,
X$ and $P_1,\ldots,P_s$ are all the singularities of~$X$.
\end{prop}

We do not prove this result fully here, but make two remarks.  Firstly, it
is clear in principle that, in the above statement, a general~$S \in \PPP$
is a quasismooth Du Val K3~surface (though in fact in the case~$\,a_2 = 1$
it is not immediate that~$S$ is quasismooth: we are allowed a
general~$X$ and a general~$S \in \PPP$ but must prove the result for every
possible~$\PPP$, not just a general choice).  Secondly, it is clear
that $C_0,\ldots,C_r \in \CS(X,\PPP)$.  This shows that
Theorem~\ref{thm:curves} excludes as many curves as it possibly could.

We say no more about how~\ref{thm:curves} is used to solve cases of
Problem~\ref{prob:main}: see~\cite{Ry06}, \cite{Ch05} and~\cite{Ry02} for
details.

\subsection*{Contents of this paper}

The remaining sections of the present paper are devoted to proving
Theorem~\ref{thm:curves}.  The proof requires several different methods and
explicit checks of dozens of cases, so often there is no choice but to give
an example calculation and a list of other cases that are similar, together
with case-specific choices that need to be made. We have therefore thought
it best to split the material up into sections according to the type of
exclusion argument used. The first of these, Section~\ref{sec:coarse},
contains arguments that are coarse and elementary --- really they are just
lemmas about curves of low degree in weighted projective 4-space --- but
they still dispose of a large number of families.  Sections~\ref{sec:tcm}
and~\ref{sec:surface} then deal with the curves that slipped through the
net, of which there are many more than one might wish.  The arguments of
Section~\ref{sec:surface} are generally more fiddly than those of
Section~\ref{sec:tcm}, and they are also required for a good many more
cases; these are summarised in Table~\ref{table:curves_surf_methods}.

\subsection*{Conventions and assumptions}

Our notations and terminology are mostly as in, for
example,~\cite{KM}, but we list here some conventions that are
nonstandard, together with assumptions that will hold throughout.

\begin{emp}
All varieties considered are complex, and they are projective and normal
unless otherwise stated.
\end{emp}

\begin{emp}
All curves are reduced and irreducible unless otherwise stated.
\end{emp}

\begin{famous95}
These are ordered as in \cite{Fl00} and \cite{CPR}, and we assume
known the basic facts about them such as quasismoothness, $\Cl X \iso
\ZZ$,~etc.  We choose coordinates $(x,y,z,t,u)$ or $(x,y_1,y_2,z,t)$
etc.\ in order of ascending degree, again as in~\cite{CPR} --- for
example, in the case of family~36, as we saw in Theorem~\ref{thm:36}, we
choose $(x_0,x_1,y,z,t)$ as coordinates for~$\PP(1,1,4,6,7)$. If~$v$~is a
coordinate then $P_v$~denotes the point where only~$v$~is nonzero.  We
import from~\cite{CPR} the notion of a \emph{starred monomial assumption\/}
--- for example, if $X = X_{15} \subset
\PP(1,1,3,4,7)_{x_0,\,x_1,\,y,\,z,\,t}$ is a member of family~25, we make
the assumption~$*tz^2$, i.e., we assume that $tz^2$~appears with nonzero
coefficient in the defining equation of~$X$. Whenever $X$~is a member of
one of the 95~families we let $A = {-K_X} = \Oh_X(1)$ denote the positive
generator of the class group; moreover, if $f \colon Y \to X$ is a
birational morphism then $B$~denotes~${-K_Y}$.
\end{famous95}

\begin{defns} \label{defns:fibr}
Let $Z$ be a normal projective variety with canonical singularities.
A~\emph{fibration\/} is a morphism $g \colon Z \to T$ to another normal
projective variety $T$ such that~$\,\dim T < \dim Z$ and~$\,g_*\Oh_Z =
\Oh_T$. We say $g$~is an \emph{elliptic fibration}, resp.\ a \emph{K3
  fibration}, if and only if its general fibre is an elliptic curve, resp.\
a K3 surface.
\end{defns}

\begin{emp}
Usually when we write an equation explicitly or semi-explicitly in
terms of coordinates we omit scalar coefficients of monomials; this is the
`coefficient convention'.
\end{emp}

\begin{emp}
If the letter~$n$ is used without explicit definition, it refers to the
degree of the mobile linear system~$\HH$ on~$X$, as in the statement of
Theorem~\ref{thm:curves}.
\end{emp}

\subsection*{Acknowledgments}

I am grateful to Miles Reid and Gavin Brown for their assistance with
several arguments that appear in this paper but date back to~\cite{Ry02},
and also to Ivan Cheltsov for some recent helpful comments.

\section{Coarse numerics and curves of low degree}
\label{sec:coarse}

Our first lemma uses the standard argument to bound the degree of a
curve centre.

\begin{lemma} \label{lemma:degC}
Let $X$ be any hypersurface in one of the 95 families and $C \subset
X$ a curve, reduced but possibly reducible.  Suppose $\HH$ is a mobile
linear system of degree~$n$ on~$X$ such that $K_X + \frac{1}{n}\HH$ is
strictly canonical and each irreducible component of~$C$ belongs to
$\CSXnH$. Then~$\,\deg C \,=\, AC \,\le\, A^3$.
\end{lemma}

\begin{proof}[\textsc{Proof}]
Let $s$ be a natural number such that $sA$ is
Cartier and very ample, and pick general members $H,H' \in \HH$. Now by
assumption
\[ \mult_{C_i}(H) \:=\: \mult_{C_i}(H') \:=\: n \]
for each irreducible component~$C_i$ of~$C$, so for a general $S \in
\left|sA \right|$
\[ A^3 sn^2 \:=\: SHH' \ge sn^2 AC \:=\: sn^2 \deg C, \]
which proves that~$\,\deg C \le A^3$.
\end{proof}

It is now necessary to understand the geometry of curves of low
degree, i.e., degree at most~$A^3$, lying inside our $X = X_d \subset
\PP(1,a_1,a_2,a_3,a_4)$. The statement of Theorem~\ref{thm:curves} suggests
the following natural case division.

\textsc{Case 1: $a_1 > 1$.} $\,\left|\Oh_X(1)\right| = \left<x_0\right>$
is fixed so there do not exist two independent degree~$1$
forms~$\,\ell,\ell'$; therefore we are trying to exclude \emph{all\/}
curves. Lemma~\ref{lemma:curvesinwps1} below shows that for many families
with $a_1 > 1$ there are in fact no curves of degree at most~$A^3$
inside~$X$, other than (perhaps) curves contracted by $\pi_4 \colon X
\dashrightarrow \PP(1,a_1,a_2,a_3)$; so for these families we have already
nearly won. There are five families with~$a_1 > 1$ to which
Lemma~\ref{lemma:curvesinwps1} does not apply, and we also need to consider
curves contracted by~$\pi_4$ --- see Lemma~\ref{lemma:curvesinwps3}, which
applies also to many families with~$a_1 = 1$, although there are
exceptional cases both with~$a_1 > 1$ and with~$a_1 = 1$ that fail to
satisfy the hypotheses.

\textsc{Case 2: $a_1 = 1$ and $a_2 > 1$}. $\,\left|\Oh_X(1)\right| =
\left<x_0,x_1\right>$ is a pencil so we are trying to exclude all
curves not contained in $\{x_0 = x_1 = 0\} \cap
X$. Lemma~\ref{lemma:curvesinwps2} below shows that for many of these
families any curve $C \subset X$ that is not contracted by~$\pi_4$ and not
contained in $\{x_0 = x_1 = 0\} \,\cap\, X$ has degree larger than~$A^3$,
so it is excluded by Lemma~\ref{lemma:degC}. Again there are families that
have $a_1 = 1$ and $a_2 > 1$ but fail to satisfy the hypothesis --- in fact
there are twelve such families --- and, as already mentioned, curves
contracted by~$\pi_4$ are considered separately in
Lemma~\ref{lemma:curvesinwps3}.

\textsc{Case 3: $a_0 = a_1 = a_2 = 1$.} Families with
$\dim\left|\Oh_X(1)\right| \ge 2$ are dealt with in the next
section.

\begin{lemma} \label{lemma:curvesinwps1}
Let $X = X_d \subset \PP = \PP(1,a_1,a_2,a_3,a_4)$ be a hypersurface
in one of the families with~$a_1 > 1$ and suppose that either
\begin{itemize}
\item[(a)] $d < a_1a_4$ or
\item[(b)] $d < a_2a_4$ and the curve $\{x = y = 0 \}\cap
X$ is irreducible (which holds for general $X$ in a family with $a_1 >
1$ by Bertini's theorem). 
\end{itemize}
Then any curve $C \subset X$ that is not
contracted by $\pi_4 \colon X \dashrightarrow \PP(1,a_1,a_2,a_3)$
has $\deg C > A^3$. Consequently $C$~is excluded absolutely by
Lemma~\ref{lemma:degC}.
\end{lemma}

\begin{rk*} \label{lemma:curvesinwps1rk} \sloppy
Out of the families with $a_1>1$, numbers~$18,\,19,\,22,\,27$ and~28 have
$d \ge a_2a_4$, so that, as written here, this lemma fails to deal with
them. (In fact we shall see in Section~\ref{sec:surface} that the
conclusion of the lemma is true for them as well.) Of the remainder, many
have $a_1a_4 \le d < a_2a_4$, which means that part~(b) of the lemma
applies to them under the generality assumption stated; this happens for
numbers~$23,\,32,\,33,\,37$, $38,\,39,\,42,\,43$, $44,\,48,\,49,\,52$,
$55,\,56,\,59,\,63$, $64,\,65,\,72,\,73,\,77$ and~89. For the rest, the
stronger form~(a) applies and no extra generality assumption is needed:
numbers~$40,\,45,\,57,\,58$, $60,\,61,\,66,\,68$, $69,\,74,\,75,\,76$,
$78,\,\ldots,\,81$, $83,\,\ldots,\,87$ and~$90,\,\ldots,\,95$.
\end{rk*}

\begin{proof}[\textsc{Proof of~\ref{lemma:curvesinwps1}}]
Most of the following proof has already appeared in~\cite{Ry06}, but we
reproduce it here for the convenience of the reader.  The part that is not
in~\cite{Ry06} is the discussion of the cases where assumption~(\ref{assn})
below fails to hold.

So suppose, contrary to the statement of the lemma, that $C \subset X$ has
$\deg C \le A^3$ and is not contracted by~$\pi_4$; let $C' \subset
\PP(1,a_1,a_2,a_3)$ be the set-theoretic image~$\pi_4(C)$. Note that $\deg
C' \le \deg C$ --- indeed, if $H$~denotes the hyperplane section of
$\PP(1,a_1,\ldots,a_4)$ and $H'$~that of $\PP(1,a_1,a_2,a_3)$, we pick
$s \ge 1$ such that $\left|sH\right|$ and $\left|sH'\right|$ are very
ample, and calculate that
\begin{eqnarray*}
s\deg C & = & (sH)C \;=\; \pi_4^*(sH')C \\
        & = & sH'(\pi_4)_*C \;=\; srH'C' \;=\; sr\deg C' \;\ge\; s\deg C',
\end{eqnarray*}
where $r \ge 1$ is the degree of the induced morphism $\pi_4|_C \colon
C \to C'$. So in fact $\deg C$ is a multiple of~$\,\deg C'$ by the
positive integer~$r$.  (The point of $\left|sH\right|$ being very ample is
that we can move it away from~$P_4$, where $\pi_4$~is undefined, and apply
the projection formula to the morphism $\pi_4|_{\PP(1,a_1,\ldots,a_4)
\smallsetminus \{P_4\}}$.)

Now we form the diagram below.
\[ \xymatrix@C=0.1cm{
 C \ar[d]    & \subset & \PP(1,a_1,a_2,a_3,a_4) \ar@{-->}[d]^{\pi_4}\\
 C' \ar[d]   & \subset & \PP(1,a_1,a_2,a_3) \ar@{-->}[d]^{\pi_3}\\
 \{*\}       & \subset & \PP(1,a_1,a_2)\\
} \]
$C'$ is contracted by $\pi_3$
 --- indeed, if its image were a curve $C''$ we would have
\[ \deg C'' \:\le\: \deg C' \:\le\: \deg C \:\le\: A^3,  \]
but~$\,A^3 = d/(a_1a_2a_3a_4) < 1/(a_1a_2)$, since~$\,d < a_3a_4\,$~in
either case (a) or (b), and on the other hand~$\,1/(a_1a_2) \le \deg C''\,$
simply because $C'' \subset \PP(1,a_1,a_2)$ --- contradiction.

For convenience we make the following assumption:
\begin{equation}
\mbox{assume that $(a_1,a_2) = 1$.} \label{assn}
\end{equation}
(We discuss at the end of the proof what to do if $(a_1,a_2) > 1$.)
(\ref{assn}) implies that the point $\{*\} \subset \PP(1,a_1,a_2)$ is,
up to coordinate change, one of
\[ \{y = z = 0\}, \quad \{y^{a_2} + z^{a_1} = x = 0\}, \quad \{x = z =
0\} \quad \mbox{and} \quad \{x = y = 0\}, \]
using the coefficient convention in~$\,y^{a_2} + z^{a_1} = 0$. It
follows that the curve $C' \subset \PP(1,a_1,a_2,a_3)$ is defined by
the same equations. In the first case, this means that $\deg C' \,=\,
1/a_3 \,>\, d/(a_1a_2a_3a_4) \,=\, A^3$, contradiction. In the second
case~$\,\deg C' \,=\, 1/a_3$ again, because
\[ C' \:\iso\: \{y^{a_2} + z^{a_1} = 0\} \:\subset\: \PP(a_1,a_2,a_3) \]
passes only through the singularity~$(0,0,1)$, using (\ref{assn}) --- so we
obtain a contradiction as in the first case. In the case $C' = \{x = z =
0\}$, we have $\deg C' = 1/(a_1a_3)$ and we easily obtain a contradiction
from $a_2a_4 > d$. In the final case,~$\,C' \,=\, \{x = y = 0\}$, if the
assumptions in part~(a) of the statement hold then we have
\[ \deg C' \:=\: 1/(a_2a_3) \:>\: d/(a_1a_2a_3a_4) \:=\: A^3, \]
contradiction; while if the assumptions in part~(b) hold then
\[ C \:=\: \{x = y = 0 \}\cap X \]
(because the right hand side is irreducible), but
\[ \deg\big(\{x = y = 0 \}\cap X \big) \:=\: a_1A^3 \:>\: A^3,  \]
since we also assumed $a_1 > 1$ --- contradiction.

This completes the proof subject to the assumption~(\ref{assn}); we
now discuss what to do if it does not hold.  First we note that
there are only 9 families with $(a_1,a_2) > 1$, namely
numbers~$18,\,22,\,28,\,43,\,52,\,59,\,69,\,73$ and~81. The first three of
these fail to satisfy either (a) or (b), so we need not concern ourselves
with them --- though we remark that the argument we are about to give works
for number~18 and fails for 22 and~28, with the inequality becoming an
equality.  Now consider as an example family~43, $X_{20} \subset
\PP(1,2,4,5,9)_{x,\,y,\,z,\,t,\,u}$ with $A^3 = 1/18$, and assume that
$\{*\} \,=\, \{y^2 + z = x = 0\} \,\subset\, \PP(1,2,4)$, which is
obviously the only problem case. Then
\[ C' \:=\: \left(\{y^2 + z = x = 0\} \,\subset\, \PP(1,2,4,5)\right)
\:\iso\: \left(\{y^2 + z = 0\} \,\subset\, \PP(2,4,5)\right), \]
which of course has
\[ \deg C' \,=\, 1/(a_3\hcf(a_1,a_2)) \,=\, 1/(5\times2) \,=\, 1/10 \,>\,
1/18, \]
contradiction. Exactly the same observation works for
numbers~$52,\,59,\,69,\,73$ and~81: one needs only to check that
$1/(a_3\hcf(a_1,a_2)) > A^3$, which is true in each case.
\end{proof}

\begin{lemma} \label{lemma:curvesinwps2}
Let $X = X_d \subset \PP = \PP(1,1,a_2,a_3,a_4)$ be a hypersurface in
one of the families with $a_1 = 1$ and $a_2 > 1$; suppose that $d <
a_2a_4$. Then any curve $C \subset X$ that is not
contracted by $\pi_4$
and that satisfies $\deg C \le A^3$ is contained in
$\{x_0 = x_1 = 0\} \cap X$.
\end{lemma}

\begin{rk*} \label{lemma:curvesinwps2rk}
Out of the families with $a_1 = 1$ and $a_2 > 1$ this lemma fails to deal
with numbers~$7,\,9,\,11,\,12,\,13,\,15,\,16,\,17,\,21,\,24,\,29$
and~34. These require extra work: see Section~\ref{sec:surface} and
particularly Table~\ref{table:curves_surf_methods}.
\end{rk*}

\begin{proof}[\textsc{Proof of~\ref{lemma:curvesinwps2}}]
Take such a curve $C$ and suppose $C \not\subset \{x_0 = x_1 = 0\}$.
\[ \xymatrix@C=0.1cm{
 C \ar[d]    & \subset & \PP(1,1,a_2,a_3,a_4) \ar@{-->}[d]^{\pi_4}\\
 C' \ar[d]   & \subset & \PP(1,1,a_2,a_3) \ar@{-->}[d]^{\pi_3}\\
 \{*\}       & \subset & \PP(1,1,a_2)\\
} \]
As in Lemma~\ref{lemma:curvesinwps1} the image of $C'$ under $\pi_3$ is a
 point --- indeed, if the image were a curve $C''$ we would have
\[ \deg C'' \:\le\: A^3 \:=\: d/(a_2a_3a_4) \:<\: 1/a_2 \:\le\: \deg C'',
 \]
because~$\,d \,<\, a_2a_4 \,\le\, a_3a_4$ --- contradiction. Therefore
after coordinate change $C' \,=\, \{x_1 = x_2 = 0\}$, since by
assumption $C' \,\ne\, \{x_0 = x_1 = 0 \}$, and so
\[ \deg C' \:=\: 1/a_3 \:>\: d/(a_2a_3a_4) \:=\: A^3,  \]
contradiction.
\end{proof}

\begin{emp} \label{wps3}
Now we need to deal with curves contracted by $\pi_4$. As discussed in
\cite[5.6]{CPR}, we can write the equation of $X$ in one of the forms
\begin{itemize}
\item[(a)] $x_4^3 + ax_4 + b \:=\: 0$, or
\item[(b)] $x_4^2        + b \:=\: 0$, or
\item[(c)] $x_jx_4^2 + ax_4 + b \:=\: 0$ (with $j = 1$, $2$ or $3$),
\end{itemize}
where $a(x_0,\ldots,x_3)$ and $b(x_0,\ldots,x_3)$ are weighted
homogeneous polynomials of the appropriate degrees. In cases (a)~and~(b),
$\pi_4 \colon X \dashto \PP(1,a_1,a_2,a_3)$ is a morphism with finite
fibres; in case~(c), $\pi_4$~contracts a finite set of curves whose union
is $\{x_j = a = b = 0 \} \,\subset\, X$.
\end{emp}

\begin{lemma} \label{lemma:curvesinwps3}
Suppose $X = X_d \subset \PP(1,a_1,\ldots,a_4)$ is a general
hypersurface in one of the 95 families and assume $d <
a_1a_2a_3$. Then any curve  $C \subset X$ contracted by $\pi_4$ has
$\deg C > A^3$, and is therefore excluded absolutely by
Lemma~\ref{lemma:degC}.
\end{lemma}

\begin{rk*}
This lemma fails to deal with families such that $P_4 \in X$ and
$\mbox{$d \ge a_1a_2a_3$}$. These are number~18, which has $a_1 > 1$;
numbers~$7,\,12,\,13$, $16,\,20,\,24$, $25,\,26$ and~46, which have $a_1 =
1$ and $a_2 > 1$; and numbers~2, 5 and~8, which have~$\,a_0 = a_1 = a_2 =
1$.
\end{rk*}

\begin{proof}[\textsc{Proof of~\ref{lemma:curvesinwps3}}]
If there exists a contracted curve $C$ then a fortiori the equation of
$X$ takes the form (c) of~\ref{wps3} above. Consider the subscheme $Z$
of the space $\PP^2(a_0,\ldots,\widehat{a_j},\ldots,a_3) =:
\PP(a_0',a_1',a_2')$ defined by $Z = \{a = b = 0\}$,  substituting
$x_j = 0$ into $a$ and $b$. $Z$ is a finite set of points (because
$a,b \in k[x_{a_0'},x_{a_1'},x_{a_2'}]$ have no common factor --- see
\cite[4.5]{CPR}) and the union of the contracted curves is the cone
over $Z$ obtained by varying $x_4$, still with $x_j = 0$. Below we
show that
\begin{equation} \label{Z} \mbox{for general
    $X$,~$\;Z$ misses any singular points of $\PP(a_0',a_1',a_2')$;}
\end{equation}
therefore our contracted curve $C$ passes through only one singular
point of $X$, namely $P_4$. Consequently
\[ \deg C \:\ge\: 1/a_4 \:>\: d/(a_1a_2a_3a_4) \:=\: A^3  \]
as required.

It remains to show (\ref{Z}).  We assume $j = 1$ to simplify the
notation --- no generality is lost in doing so because the proof below
does not make use of~$\,a_1 \le a_2 \le a_3$.  We know
\[ Z \:=\: \{a_{d-a_4} = b_d = 0 \} \:\subset\: \PP(1,a_2,a_3)  \]
and we need to show that either~$\,a_2 \newmid (d - a_4)$ or~$\,a_2 \newmid
d$.  This demonstrates that $(0,1,0) \notin Z$, assuming $X$ is
general. Formally we also need to show that either~$\,a_3 \newmid (d -
a_4)$ or~$\,a_3 \newmid d$, but the proof is identical. Note that even
if~$(a_2,a_3) \ne 1$ the only two points of $\PP(1,a_2,a_3)$ which can be
singular are~$(0,1,0)$ and~$(0,0,1)$.

Now to the proof. Because $x_1x_4^2$ is the tangent monomial to $X$ at
$P_4$, we know that
\begin{eqnarray}
a_1 + 2a_4 & = & d \qquad \mbox{and} \label{Z1}\\
a_2 + a_3 & = & a_4, \label{Z2}
\end{eqnarray}
where (\ref{Z2}) follows from (\ref{Z1}) and~$\,d \,=\, a_1 + \cdots +
a_4$. Now we consider the different possibilities for the tangent
monomial to $X$ at~$P_2$.

If~$\,x_4x_2^n$ is the tangent monomial to~$X$ at~$P_2$ then~$\,a_4 + na_2
\,=\,d$, so $\mbox{$a_2 \newmid (d-a_4)$}$ and we are done. If
$x_3x_2^n$~is the tangent monomial at $P_2$ then $a_3 + na_2 \,=\, d$, so
\[ (n-1)\,a_2 \:=\: d-a_4 \] using~(\ref{Z2}), which shows that~$\,n \ge 2$
and~$\,a_2 \newmid (d-a_4)$ as required. If $x_2^n$ is the tangent monomial
then $P_2 \notin X$ and~$\,a_2 \newmid d$.

We are left with the case $x_1x_2^n$. We know that
\begin{eqnarray}
a_1 + na_2 & = & d \qquad \qquad \mbox{and} \label{Z3}\\
a_3 + a_4 & = & (n-1)\,a_2, \label{Z4}
\end{eqnarray}
where as before (\ref{Z4})~follows from~(\ref{Z3}) and~$\,d \:=\: a_1 +
\cdots + a_4$. Now (\ref{Z2})~and~(\ref{Z4}) imply
\begin{eqnarray*}
2a_3 & = & (n-2)\,a_2 \qquad \mbox{and}\\
2a_4 & = & na_2.
\end{eqnarray*}
If $n$ is even then~$\,a_2 \newmid a_3$ and~$\,a_2 \newmid a_4$, so~$\,a_2
= 1$ (any three of $(a_1,a_2,a_3,a_4)$ have highest common factor 1 because
the K3 section $\{x_0=0\} \cap X$ is well formed). Therefore~$\,a_2 \newmid
d$, as required. If on the other hand $n$ is odd then $a_2 = 2a_2'$ is even
and $a_2'$ divides $a_2$, $a_3$ and $a_4$, so $a_2' = 1$ and
\[ (a_0,\ldots,a_4) \:=\: (1,a_1,2,a_4-2,a_4)  \]
with~$\,a_4 = n$ odd. If $a_1$ is even then, by (\ref{Z3}), $d$ is even
and~$\,a_2 = 2 \newmid d$; but if $a_1$ is odd then $d$ is odd as well
and~$\,a_2 \,=\, 2 \newmid (d-a_4)$.
\end{proof}

\section{The test class method} \label{sec:tcm}

The following lemma is completely general and elementary; we will use it
for curves inside~$X$, but it is also important for excluding singular
points: see~\cite[Theorem~3.20]{Ry06}.  It should be compared
with~\cite[5.2.1]{CPR}, to which it is closely related.

\begin{lemma} \label{lemma:testcl_1}
Let $X$ be a Fano 3-fold hypersurface in one of the 95~families and
$\HH$ a mobile linear system of degree~$n$ on~$X$ such that $K_X +
\frac{1}{n}\HH$ is strictly canonical; suppose $\Gamma\subset X$ is an
irreducible curve or a closed point satisfying $\Gamma\in\CSXnH$, and
furthermore that there is a Mori extremal divisorial contraction
\[ f\colon (E\subset Y) \to (\Gamma\subset X), \quad \Centre_X E = \Gamma,
\]
such that $E\in\VXnH$ (for the notation $\VXnH$, see the Introduction).
Then $B^2 \in \NEbar Y$.
\end{lemma}

\begin{proof}[\textsc{Proof}]
We know that
\[ K_Y + \textstyle\frac{1}{n}\HH_Y \:\qeq\: f^*\left(K_X +
\textstyle\frac{1}{n}\HH\right) \:\qeq\: 0 . \]
It follows that $B \qeq \frac{1}{n}\HH_Y$, and therefore $B^2 \in
\NEbar Y$, because $\HH_Y$ is~mobile.
\end{proof}

The idea of the test class method is very simple. Suppose $\Gamma
\subset X$ is an irreducible curve or a closed point that is the
centre of a Mori extremal divisorial contraction
$f \colon (E \subset Y) \to (\Gamma \subset X)$
as in the above lemma. A~\emph{test class\/} is, by definition, a
nonzero nef class~$M \in \None\! Y$.

\begin{lemma}[cf. {\cite[5.2.3]{CPR}}] \label{lemma:testcl_2}
Suppose that, in the situation just described, there is a test class~$M$
on~$Y$ with~$MB^2 < 0$. Then $E$~cannot be a strictly canonical singularity
for any~$\HH$.
\end{lemma}

\begin{proof}[\textsc{Proof}]
This is immediate from Lemma~\ref{lemma:testcl_1}.
\end{proof}

\begin{cor} \label{curvetccor}
If the hypotheses of Lemma~\ref{lemma:testcl_2} are satisfied by some curve
$C = \Gamma \subset X$ then $C$ is excluded absolutely, that is, $C$~is not
a strictly canonical centre for any~$\HH$.
\end{cor}

\begin{proof}[\textsc{Proof}]
We are assuming that there exists a Mori extremal divisorial contraction
$f \colon (E \subset Y) \to (C \subset X)$ with $\Centre_X(E) = C$.
Suppose $\HH$~is mobile of degree~$n$ on~$X$ with $K_X + \frac{1}{n}\HH$
strictly canonical. Clearly what we need to prove is the following: if $C
\in \CSXnH$ then in fact $E \in \VXnH$.  To see this, first note that over
a general point $P \in C \subset X$, $f \colon Y \to X$ must be the blowup
of~$\II_C$. Let $P \in S \subset X$ be a general surface through~$P$,
smooth near~$P$ and transverse to~$C$. Then
\[ \mult_P\left(\HH|_S\right) = n \]
because $C \in \CSXnH$ by assumption and we have the classical fact
that, locally over $P \,=\, C \cap S \,\subset\, S$, the first ordinary
blowup extracts a divisor of maximal multiplicity for~$\HH|_S$.
\end{proof}

The problem with the test class method is that it only applies to
curves $C \subset X$ that are centres of Mori extremal divisorial
contractions. Such curves are always contained in~$\,\Nonsing(X)$
and their own singularities are also restricted. It turns out that the
test class method, together with coarse arguments like those of
Section~\ref{sec:coarse}, is sufficient to prove
Theorem~\ref{thm:curves} for families with $a_0 = a_1 = a_2 = 1$; for the
other families, the curves that the coarse results fail to deal with hit
singularities of~$X$, and we need other methods.

We now turn to the more practical question of how to find a test class for
a given curve.

\begin{defn}[cf. {\cite[5.2.4]{CPR}}] \label{Gammaisol}
Let $L$ be a Weil divisor class in a 3-fold~$X$ and $\Gamma \subset X$
an irreducible curve or a closed point. We say that \emph{$L$
  isolates~$\Gamma$}, or is a \emph{$\Gamma$-isolating class\/}, if and
only if there exists $s \in \ZZ_{\ge1}$ such that the linear system
$\mbox{$\LL_{\Gamma}^s := \left|\II_{\Gamma}^s(sL) \right|\,$}$~satisfies
\begin{itemize}
\item $\Gamma \in \Bs \LL_{\Gamma}^s$ is an isolated component (i.e.,
in some neighbourhood of~$\Gamma$ the subscheme $\Bs \LL_{\Gamma}^s$
is supported on~$\Gamma$); and
\item if $\Gamma$ is a curve, the generic point of~$\Gamma$ appears
with multiplicity~1 in~$\,\Bs \LL_{\Gamma}^s$.
\end{itemize}
\end{defn}

\begin{lemma} \label{lemma:if_isol_then_test}
Suppose that $L$ isolates $\Gamma \subset X$ and let~$\,s \in \ZZ_{\ge
  1}\,$~be as above. Then for any extremal divisorial contraction
\[ f \colon (E \subset Y) \to (\Gamma \subset X) \mbox{$\:$ with $\:$}
  \Centre_X(E) = \Gamma \]
the birational transform $M = f^{-1}_*\LL^s_{\Gamma}$ is a test class
on~$Y$.
\end{lemma}

\begin{proof}[\textsc{Proof}]
This is~\cite[5.2.5]{CPR}.
\end{proof}

We now use the test class method, together with some elementary
arguments in the style of Lemmas~\ref{lemma:curvesinwps1}
and~\ref{lemma:curvesinwps2}, to prove Theorem~\ref{thm:curves} for all the
families with $a_0 = a_1 = a_2 = 1$, that is, for families
$1,\ldots,6,8,10$ and~$14$.

\begin{proofofcthm31s}
Let
\[ X \:=\: X_d \:\subset\: \PP(1,1,1,a_3,a_4) \]
be a hypersurface in one of families~$1,\ldots,6,8,10$ and~14 and $C
\subset X$ a curve; suppose that $C$~is a strictly canonical centre
for some~$\HH$. By Lemma~\ref{lemma:degC}, $\deg C \,\le\, A^3$.

\textsc{Case 1: $C$ is contracted by $\pi_4 \colon X
\dashrightarrow \PP(1,1,1,a_3)$.} By Lemma~\ref{lemma:curvesinwps3}, we are
in a family with $d \ge a_1a_2a_3$ and $P_4 \in X$, that is, one of
families~2, 5 and~8. It is very easy to check in each of these cases that
the contracted curves are contained in $\{\ell = \ell' = 0\} \cap X$ for
two linearly independent forms~$\ell,\ell'$ of degree~$1$
in~$(x_0,\ldots,x_4)$ --- for example, in the case of family~$8$,
$\mbox{$X_9 \subset \PP(1,1,1,3,4)_{x_0,\,x_1,\,x_2,\,y,\,z}$}$ with $A^3 =
3/4$, we do a coordinate change so that the tangent monomial at $P_4 = P_z$
is~$\,x_2 z^2$; then the equation of~$X$ is
\[ x_2z^2 + a_5z + b_9 \,=\, 0 \quad \mbox{with} \quad a,b \in
k[x_0,x_1,x_2,y] \]
and the contracted curves are the irreducible components of
\[ \{x_2 = a_5 = b_9 = 0\} \:\subset\: \PP(1,1,1,3). \]
But~$\,y^3 \in b_9\,$~by quasismoothness at~$P_y$ and therefore after a
coordinate change
\[ C \:=\: \{x_1 = x_2 = y = 0\} \:\subset\: \{x_1 = x_2 = 0\} \,\cap\,
X. \]

\textsc{Case 2: $C$ is not contracted by $\pi_4$.}  As in the proofs of
Lemmas~\ref{lemma:curvesinwps1} and~\ref{lemma:curvesinwps2} we consider the
following diagram.
\[ \xymatrix@C=0.1cm{
 C \ar[d]    & \subset & \PP(1,1,1,a_3,a_4) \ar@{-->}[d]^{\pi_4}\\
 C' \ar[d]   & \subset & \PP(1,1,1,a_3) \ar@{-->}[d]^{\pi_3}\\
 C''         & \subset & \PP(1,1,1)\\
} \]
We may assume that $C'$ is not contracted by~$\pi_3$ --- indeed, any
point in~$\PP^2$ is defined by two linearly independent forms~$\ell,\ell'$
 of degree~$1$ in~$(x_0,x_1,x_2)$. So~$\,\deg C'' \in \ZZ_{\ge 1}$ and
 therefore~$\,\deg C',\,\deg C'' \in \ZZ_{\ge 1}$ also (they are positive
 integral multiples of $\deg C''$ --- see the proof of
 Lemma~\ref{lemma:curvesinwps1}). For families~8, 10~and~14, $A^3 < 1$ and
 we already have our contradiction; families~$1,\ldots,6$ remain.

The next step is to show that if $C$ is not contained in some $\{\ell
= \ell' = 0\}$ then, after a coordinate change, it is one of the
following; here $N$~denotes the number of the family.
\begin{itemize}
\item[1.] $\,N = 1$, $\,X_4 \subset \PP^4$, $\,A^3 = 4$, $\,C =$ a twisted
  cubic curve in some linearly embedded $\PP^3 \subset \PP^4$, test
  class~$2A - E$.
\item[2.] $\,N = 2$, $\,X_5 \subset \PP(1,1,1,1,2)$, $\,A^3 = 5/2$,
$\,C = \{y = x_3 = x_0x_1 + x_2^2 = 0\}$, $\,\deg C = 2$, test class~$2A -
E$.
\item[3.] $\,N = 3$, $\,X_6 \subset \PP(1,1,1,1,3)$, $\,A^3 = 2$,
$\,C = \{y = x_3 = x_0x_1 + x_2^2 = 0\}$, $\,\deg C = 2$, test class~$6A -
E$.
\item[4.] $\,N = 4$, $\,X_6 \subset \PP(1,1,1,2,2)$, $\,A^3 = 3/2$,
$\,C = \{y_2 = y_1 = x_0 = 0\}$, $\,\deg C = 1$, test class~$2A - E$.
\item[5.] $\,N = 5$, $\,X_7 \subset \PP(1,1,1,2,3)$, $\,A^3 = 7/6$,
$\,C = \{z = y = x_0 = 0\}$, $\,\deg C = 1$, test class~$6A - E$.
\item[6.] $\,N = 6$, $\,X_8 \subset \PP(1,1,1,2,4)$, $\,A^3 = 1$,
$\,C = \{z = y = x_0 = 0\}$, $\,\deg C = 1$, test class~$4A - E$.
\end{itemize}
As an illustration of how to derive this list we consider family~4 (the
others being easier). We have $X = X_6 \subset
\PP(1,1,1,2,2)_{x_0,\,x_1,\,x_2,\,y_1,\,y_2}$ with $A^3 = 3/2$;
if necessary we do a coordinate change so that $P_4 = P_{y_2} \not\in
X$. Take a curve $C \subset X$ of degree at most $A^3 = 3/2$ in~$\PP^2$, so
it is a line $\{x_0 = 0\}$ after coordinate change and
\[ \deg C = \deg C' = \deg C'' = 1. \]
Now $C' \,\subset\, \big(\{x_0 = 0\} \cap \PP(1,1,1,2)\big) \,\iso\,
\PP(1,1,2)\,$ is an irreducible curve so, after coordinate change, it is
$\{y_1 = x_0 = 0\}$. But
\[ \{y_1 = x_0 = 0\} \cap X \;\iso\; \{y_2^3 + a_2y_2^2 + b_4y_2 + c_6 \,=\,
0\} \;\subset\; \PP(1,1,2) \]
with $a,b,c \in k[x_1,x_2]$. Because $C$~is a degree~$1$ component of
this, it corresponds to a linear factor of the cubic in~$y_2$, so
after another coordinate change $C \,=\, \{y_2 = y_1 = x_0 = 0\}$, as
required.

Finally the curves in the above list need to be excluded using the test
class method. The method is essentially the same in each case, so we give
the details only for case~3. So let
$X = X_6 \subset \PP(1,1,1,1,3)_{x_0,\ldots,\,x_3,\,y}$ be a general
member of family~3 and suppose $C \,=\, \{y = x_3 = x_0x_1 + x_2^2 = 0\}
\subset X$. It is clear that $6A$~is $C$-isolating
(Definition~\ref{Gammaisol}, using~$s = 1$), so by
Lemma~\ref{lemma:if_isol_then_test} $M = 6A - E$ is a test class, where $f
\colon (E \subset Y) \to (C \subset X)$ is the blowup of~$C$. But
\begin{eqnarray*}
MB^2 & = & (6A - E)(A - E)(A - E) \\
     & = & 6A^3 \,-\, 13A^2E \,+\, 8AE^2 \,-\, E^3 \\
     & = & 6\times2 \,-\, 0 \,-\, 8\times2 \,-\, 0 \;=\; {-4} \;<\; 0
\end{eqnarray*}
so $C$ is excluded by Corollary~\ref{curvetccor}. In the calculation
we used
\[ A^2E \;=\; 0, \quad AE^2 \;=\; {-\deg C} \;=\; {-2}, \]
\[ E^3 \;=\; {-\deg\NNN_{C\mid X}} \;=\; {-\deg C} \,+\, 2 \,-\, 2p_a(C)
\;=\; 0. \]
This completes the proof.
\hfill $\qedsymbol$
\end{proofofcthm31s}

\section{Surface methods for the remaining curves}
\label{sec:surface}

The task that remains is to prove Theorem~\ref{thm:curves} for families
with~$\:a_2 > 1$.  This involves checking lots of cases; before listing
them we consider two families in full detail so as to illustrate the two
main methods we need.

\begin{ex}[Theorem~\ref{thm:curves} for family~20] \label{ex:curves_surf_1}
Take a general
\[ X_{13} \:\subset\: \PP\,(1,1,3,4,5)_{x_0,\,x_1,\,y,\,z,\,t} \mbox{\quad
  with $A^3 = 13/60$;} \]
we make two starred monomial assumptions: $*tz^2$~and~$*zy^3$. The presence
of~$yt^2$, on the other hand, is guaranteed by quasismoothness at~$P_t$, so
after a coordinate change we can write the defining equation of~$X$ as
\[ yt^2 + a_8t + b_{13} \:=\: 0 \quad \mbox{with} \quad a,\,b \in
k[x_0,x_1,y,z]. \]
Let $C$ be one of the curves of degree~$\,1/5 < 13/60 = A^3$ which
pass through~$P_t$ and are flopped by the quadratic involution~$i_{P_t}$:
so $C$~is a component of $\{y = a_8 = b_{13} = 0\} \,\subset\, X$.  The
components of this set are the only curves that remain to be excluded for
this~$X$ --- indeed,~$\,d < a_2\,a_4\,$~so Lemma~\ref{lemma:curvesinwps2}
applies to curves not contracted by~$\pi_4$ --- and of course all the
components of $\{y = a_8 = b_{13} = 0\}$ are the same up to coordinate
change, so it is enough to exclude one of them, our~$C$.

After a coordinate change we may assume $C \,=\, \{x_0 = y = z = 0\}
\,\subset\, X$, that is, $C$~is the~$\,x_1t\,$-stratum (note that~$\,z^2
\in a_8\,$~by $*tz^2$, so~$\,P_z \not\in C$ before the change). To
exclude~$C$ we follow the general method described in~\cite[\S5]{CM},
taking a general surface $T \in \left|4A - C\right|$ and doing the
following calculations.

\begin{claim*} (a) $\Bs\left|4A - C\right|$ is supported on~$\,C \cup
  \{P_y\}$.

(b) $T$ has a $\frac{1}{5}(1,1)$ singularity at~$\,P_t \in C \subset
T$ and $T$~is smooth at all other points of~$C$.

(c) The selfintersection~$(C)^2_T = {-9/5}$.
\end{claim*}

\noindent For the proof, see below. Suppose now that $\HH$~is a mobile
linear system of degree~$n$ on~$X$ such that $K_X + \frac{1}{n}\HH$ is
strictly canonical and $C \in \CSXnH$. Restricting to~$T$, we have $\HH|_T
= nC + \LL$ where $\LL$~is the mobile part. It follows that
$\left(\frac{1}{n}\HH|_T - C\right) \,\qeq\, \frac{1}{n}\LL\,$~is nef
on~$T$; but we~calculate
\begin{eqnarray*}
\left(\textstyle\frac{1}{n}\LL\right)^2_T \;=\;
\left(\textstyle\frac{1}{n}\HH|_T - C\right)^2_T & = &
\left(A|_T\right)^2 \,-\, 2\left(A|_T\right)C \,+\, (C)^2_T \\
 & = & A^2T \,-\, 2AC \,+\, (C)^2_T \\
 & = & 4\times\textstyle\frac{13}{60} \,-\,
 2\times\textstyle\frac{1}{5} \,-\, \textstyle\frac{9}{5} \;=\;
 {-\textstyle\frac{4}{3}} \;<\; 0,
\end{eqnarray*}
contradiction.

\begin{proof}[\textsc{Proof of claim}]
(a) A general element $T \in \left|4A - C\right|$ has equation
\begin{equation}
z \,+\, yS^1(x_0,x_1) \,+\, x_0\,S^3(x_0,x_1) \;=\; 0, \label{eqnT}
\end{equation}
with the coefficient convention. If $P \in \Bs\left|4A - C\right|$
then clearly $z = x_0 = 0$ at $P$; if $y \ne 0$ then $x_1 = 0$ so~$\,a_8 =
b_{13} = 0$ because neither contains a pure power of~$y$, and it
follows from the defining equation of~$X$ that~$\,t=0$.

(b) Inside $X$, $P_t \sim \frac{1}{5}(1,1,4)$ in local coordinates
$(x_0,x_1,z)$. The usual manipulation of the defining equation of~$X$,
together with a local analytic coordinate change, shows that
$y \,=\, z^2 + x_0^8 + \cdots + x_0\,x_1^7\,$
near~$P_t$ (note that $x_1^8$~does not appear because~$\,C \subset
X$). Therefore a general $T \in \left|4A - C\right|$, which is
globally defined by~(\ref{eqnT}), is locally defined by
\[ z \:+\: \left(z^2 + x_0\,S^7(x_0,x_1)\right)S^1(x_0,x_1) \:+\:
x_0\,S^3(x_0,x_1) \;=\; 0, \]
so $(P_t \in T) \sim \frac{1}{5}(1,1)$ in local
coordinates~$(x_0,x_1)$. Note that near~$P_t \in T$ the curve~$C$ is
defined by~$\,x_0 = 0$.

To show that $T$ is smooth at all other points of $C$ one considers
the affine piece $\{x_1 \ne 0\} \subset \PP(1,1,3,4,5)$, inside which
$T$ is defined by
\[ yt^2 + a_8t + b_{13} \;=\; 0 \quad \mbox{and} \quad z + y + yx_0 + x_0^4
+ \cdots + x_0 \;=\; 0 \]
with $a,b \in k[x_0,y,z]$. Writing down the four partial derivatives
of each of these two expressions, and evaluating them along~$\,\{x_0 = y
= z = 0\}$, one sees that if $X$~is general the rank of the $4 \times
2$~matrix never drops below~$2$.

(c) The nontrivial part here is to calculate the \emph{different},
$\Diff \subset C$, which is the divisor satisfying
\[ \left(K_T + C\right)|_C \;=\; K_C + \Diff\! . \]
$C$ is Cartier away from $P_t \in T$ so~$\,\Diff$ is supported on~$P_t$
and the only problem is to calculate the coefficient. We use Corti's
result~\cite[16.6.3]{FA}, which implies that $\Diff = \frac{m-1}{m}P_t$
where $m$~is the index of~$C$ at $P_t \in T$, provided that $K_T + C$
is plt (purely log terminal) at~$P_t$. But the plt condition is clear in
this case: $P_t \in T$ is resolved by the $\frac{1}{5}(1,1)$ (i.e.,
ordinary) blowup, the discrepancy of~$K_T$ is $1/5 - 4/5 = {-3/5}$ (because
$a_E(K_X) = 1/5$ for the $\frac{1}{5}(1,1,4)$ blowup of~$P_t \in X$, and
$T$~has local weight $4/5$), and $C \subset T$ has local weight~$1/5$; so
the log discrepancy of $K_T + C$ is ${-3/5} - 1/5 \,=\, {-4/5} \,>\,
{-1}$. Clearly~$\,m = 5$, so~$\,\Diff = \frac{4}{5}P_t$.

The rest is easy. $T \subset X$ is Cartier in codimension~$2$, because
$X$~has isolated singularities, so
\begin{eqnarray*}
{-2} \,+\, \textstyle\frac{4}{5} & = & \left(K_T \,+\, C\right)C \;=\;
(C)^2_T \,+\, \left(K_X \,+\, T\right)C \\
 & = & (C)^2_T \,+\, 3AC \;=\; (C)^2_T \,+\, \textstyle\frac{3}{5}
\end{eqnarray*}
and therefore $(C)^2_T = {-9/5}$ as required.
\end{proof}
\end{ex}

\begin{ex} \label{ex:curves_surf_2}
Family~29, $X_{16} \subset \PP(1,1,2,5,8)_{x_0,\,x_1,\,y,\,z,\,t}$
with~$\,A^3 = 1/5$. Suppose that $X$~contains the curve $C \,=\, \{x_0 = y
= t = 0\}$. (An easy argument in the style of the proofs of
Lemmas~\ref{lemma:curvesinwps1} and~\ref{lemma:curvesinwps2} shows that up
to coordinate change this $C$ is the only curve of degree at most~$A^3$
not contained in $\{x_0 = x_1 = 0\}$.) We can write the equation of~$X$ as
\[ t^2 + a_8t + b_{16} \;=\; 0 \quad \mbox{with} \quad a,\,b \in
k[x_0,x_1,y,z]. \]
We have assumed that $C \subset X$, which means that after
making the substitution $x_0 = y = 0$ in $a$~and~$b$ we are left with a
reducible quadratic $t(t + c_8) = 0$, where $c \in k[x_1,z]$. In other
words,
\[ \Bs|2A - C| \;=\; \{x_0 = y = 0\} \cap X \;=\; C + C', \]
where $C' \,=\, \{x_0 = y = t + c_8 = 0\}$ is just like $C$ after a
coordinate change. Now let $T \in |2A - C|$ be a general surface.

\begin{claim*}
(a) $T$ has a $\frac{1}{5}(1,3)$ singularity at $P_z$ and is
smooth elsewhere.

(b) The selfintersection $(C)^2_T = {-7/5}$ and so, by
symmetry, $(C')^2_T = {-7/5}$ as well.
\end{claim*}

See below for the proof. Suppose now that $\HH$~is mobile
of degree~$n$ on~$X$ with $\KXnH$~canonical and $C \in \CSXnH$. Then,
restricting to~$T$, we have $\frac{1}{n}\HH|_T \,\qeq\, C + \alpha C' +
\frac{1}{n}\LL$, where $0 \le \alpha \le 1$ and $\LL$~is the mobile part
of~$\HH|_T$. But $\frac{1}{n}\HH|_T \,\qeq\, A|_T \,=\, C + C'$,
so~$\,\frac{1}{n}\LL \,\qeq\, (1 - \alpha)C'$. It follows that
\[ 0 \:\le\: (\textstyle\frac{1}{n}\LL)^2_T \:=\: (1 - \alpha)^2(C')^2_T
\:=\: {-\frac{7}{5}}(1 - \alpha)^2, \]
and therefore~$\alpha = 1$.

Consequently $C' \in \CSXnH$ as well --- but~$\,\deg(C + C') = 2A^3$,
contradicting Lemma~\ref{lemma:degC}.

\begin{proof}[\textsc{Proof of Claim}]
(a) Near to $P_z$, after a local analytic coordinate change, $T \,=\, \{y
= 0\}$, so clearly $P_z \sim \frac{1}{5}(1,3)$ inside~$T$. Showing $T$~is
smooth elsewhere can be done as in the proof of the claim in
Example~\ref{ex:curves_surf_1} above.

(b) This is also essentially the same as the calculation in
Example~\ref{ex:curves_surf_1}.  We check that $K_T + C$~is plt at~$P_z$, and
it is clear that the index of~$C$ at~$P_z$ is~5, so by Corti's
result~\cite[16.6.3]{FA} we have
\begin{eqnarray*} {-2} + \textstyle\frac{4}{5} \;=\; \deg(K_C + \Diff) & =
  & (K_T + C)C \;=\; (K_X + T)C + (C)^2_T \\
& = & AC + (C)^2_T \;=\; \textstyle\frac{1}{5} + (C)^2_T.
\end{eqnarray*}
The desired conclusion follows.
\end{proof}

\end{ex}

\begin{proofofcthm11}
For the majority of the families with $a_1 > 1$, Lemmas~\ref{lemma:degC},
\ref{lemma:curvesinwps1} and~\ref{lemma:curvesinwps3} prove
Theorem~\ref{thm:curves}. We need consider only
families~$18,\,19,\,22,\,27$ and~28, which fail to satisfy the hypotheses
for Lemma~\ref{lemma:curvesinwps1} --- and family~18 also fails to satisfy
the hypotheses for Lemma~\ref{lemma:curvesinwps3}. The way things turn out
is as follows: firstly, for families~$19,\,22,\,27$ and~28 there are in
fact no curves of degree at most~$A^3$ contained in~$X$, and for
family~18 the only curves of degree at most~$A^3$ are those contracted
by~$\pi_4$ --- in other words, Lemma~\ref{lemma:curvesinwps1} in fact
applies to \emph{all\/} the families with $a_1 > 1$, provided we make
generality assumptions. Secondly, the curves in family~18 contracted
by~$\pi_4$ can be excluded as in Example~\ref{ex:curves_surf_1}, using a
general surface~$T \in |4A - C|$.

We make no further remarks about the exclusion of curves contracted
by~$\pi_4$ in the case of family~18, but give an example of how to extend
Lemma~\ref{lemma:curvesinwps1} to families~$18$, 19, $22,\,27$ and~28. So
consider family~19,
\[ \mbox{$X_{12} \:\subset\: \PP(1,2,3,3,4)_{x,\,y,\,z_1,\,z_2,\,t}\;$ with
  $\;A^3 = 1/6$.} \]
Let $P_1,P_2,P_3,P_4 \sim \frac{1}{3}(1,2,1)$ be the
singularities on the $\mbox{$z_1z_2$-stratum}$ and $Q_1,Q_2,Q_3 \sim
\frac{1}{2}(1,1,1)$ those on the $yt$-stratum. We assume that the
curve $\{x = y = 0\} \cap X$ is irreducible and that $P_iQ_j
\not\subset X$ for all~$i,\,j$; a general~$X$ satisfies these
assumptions. Now suppose that $C$ is a curve of degree at most~$A^3$
contained in~$X$. Again we form the following familiar diagram.
\[ \xymatrix@C=0.1cm{
 C \ar[d]    & \subset & \PP(1,2,3,3,4) \ar@{-->}[d]^{\pi_4}\\
 C' \ar[d]   & \subset & \PP(1,2,3,3) \ar@{-->}[d]^{\pi_3}\\
 C''         & \subset & \PP(1,2,3)\\
} \]
Certainly $C'$ is a curve, because $P_t \not\in X$. Suppose that $C''$~is
also a curve. Then its degree is~$1/6 = A^3$ and it is defined by $\{x =
0\}$ after a coordinate change. Therefore $\deg C = \deg C' = 1/6$ also,
and $C'$~is isomorphic to a curve in~$\PP(2,3,3)$, so after a coordinate
 change we have $C' \,=\, \{x = z_1 = 0\}$. The same argument applied
 to~$C$ now shows that $C \,=\, \{x = z_1 = t = 0\}$, after another
 coordinate change, so $C = P_iQ_j$,~contradiction.

Therefore in fact $C'' = \{*\}$ is a point. After a coordinate change,
this point is one of
\[ \{y = z_1 = 0\}, \quad \{x = y^3 + z_1^2 = 0\}, \quad
\{x = z_1 = 0\} \quad \mbox{and} \quad \{x = y = 0\}. \]
In the first two cases $\deg C' = 1/3 > A^3$, contradiction. In the
last case, $C \,\subset\, \{x = y = 0\} \cap X$, but the right-hand side
has degree~$1/3 > A^3$ and is irreducible by assumption ---
contradiction again. In the third case, an easy argument shows that $C
= P_iQ_j$ for some~$i,\,j$, which gives a contradiction as above.

Similar arguments can be used to extend Lemma~\ref{lemma:curvesinwps1} to
families~$18,\,22$, 27~and~28. This completes the proof of
Theorem~\ref{thm:curves} for all the families with~$\mbox{$a_1 >
  1$}$. \hfill $\qedsymbol$
\end{proofofcthm11}

\begin{proofofcthm21s} \label{prfofcurvesthm21s}
For this proof we apply the method of Example~\ref{ex:curves_surf_1} to many
curves.  There is not enough space here to go through each of these;
instead, Table~\ref{table:curves_surf_methods} summarises the calculations.
\begin{table}[htb] \caption{Curves excluded by surface methods}
\label{table:curves_surf_methods}
\begin{center}
\begin{tabular}{|l|l|l|l|l|}
\hline
Family & Fails & Curve(s) & Method & System \\
\hline
7 & \ref{lemma:curvesinwps2}, \ref{lemma:curvesinwps3} & $\{x_0 = y_1 = y_2 = 0\}$
& \ref{ex:curves_surf_1} & $|2A - C|$ \\
 & & $\{x_0 = y_1 = z = 0\}$ & \ref{ex:curves_surf_1} & $|3A - C|$ \\
\hline
9 & \ref{lemma:curvesinwps2} & $\{x_0 = y = z_1 = 0\}$
& \ref{ex:curves_surf_1} & $|3A - C|$ \\
 & & $\{x_0 = z_1 = z_2 = 0\}$ & \ref{ex:curves_surf_1} & $|3A - C|$ \\
\hline
11 & \ref{lemma:curvesinwps2} & $\{x_0 = y_1 = z = 0\}$
& \ref{ex:curves_surf_1} & $|5A - C|$ \\
\hline
12 & \ref{lemma:curvesinwps2}, \ref{lemma:curvesinwps3} & $\{x_0 = y = z = 0\}$
& \ref{ex:curves_surf_1} & $|3A - C|$ \\
 & & $\{x_0 = y = t = 0\}$ & \ref{ex:curves_surf_1} & $|4A - C|$ \\
\hline
13 & \ref{lemma:curvesinwps2}, \ref{lemma:curvesinwps3} & $\{x_0 = y = z = 0\}$
& \ref{ex:curves_surf_1} & $|3A - C|$ \\
 & & $\{x_0 = y = t = 0\}$ & \ref{ex:curves_surf_1} & $|5A - C|$ \\
\hline
15 & \ref{lemma:curvesinwps2} & $\{x_0 = y = t = 0\}$
& \ref{ex:curves_surf_2} & $|2A - C|$ \\
\hline
16 & \ref{lemma:curvesinwps2}, \ref{lemma:curvesinwps3} & $\{x_0 = y = z = 0\}$
& \ref{ex:curves_surf_1} & $|4A - C|$ \\
\hline
17 & \ref{lemma:curvesinwps2} & $\{x_0 = y = z_1 = 0\}$
& \ref{ex:curves_surf_1} & $|4A - C|$ \\
\hline
20 & \ref{lemma:curvesinwps3} & $\{x_0 = y = z = 0\}$
& \ref{ex:curves_surf_1} & $|4A - C|$ \\
\hline
21 & \ref{lemma:curvesinwps2} & $\{x_0 = y = t = 0\}$
& \ref{ex:curves_surf_1} & $|7A - C|$ \\
\hline
24 & \ref{lemma:curvesinwps2}, \ref{lemma:curvesinwps3} & $\{x_1 = y = z = 0\}$
& \ref{ex:curves_surf_1} & $|5A - C|$ \\
\hline
25 & \ref{lemma:curvesinwps3} & $\{x_1 = y = z = 0\}$
& \ref{ex:curves_surf_1} & $|4A - C|$ \\
\hline
26 & \ref{lemma:curvesinwps3} & $\{x_0 = y = z = 0\}$
& \ref{ex:curves_surf_1} & $|5A - C|$ \\
\hline
29 & \ref{lemma:curvesinwps2} & $\{x_0 = y = t = 0\}$
& \ref{ex:curves_surf_2} & $|2A - C|$ \\
\hline
34 & \ref{lemma:curvesinwps2} & $\{x_0 = y = t = 0\}$
& \ref{ex:curves_surf_2} & $|2A - C|$ \\
\hline
46 & \ref{lemma:curvesinwps3} & $\{x_1 = y = z = 0\}$
& \ref{ex:curves_surf_1} & $|7A - C|$ \\
\hline
\end{tabular}
\end{center}
\end{table}

The contents of Table~\ref{table:curves_surf_methods} should be interpreted
as follows, in conjunction with the Big Table of~\cite{CPR}. The families
listed are those with $a_1 = 1$ and $a_2 > 1$ which fail to satisfy the
hypotheses of at least one of Lemmas~\ref{lemma:curvesinwps2}
and~\ref{lemma:curvesinwps3}; which of these two they fail is the content
of the second column. Now, for a given family in the table, we run familiar
arguments, in the style of the proofs of Lemmas~\ref{lemma:curvesinwps1}
and~\ref{lemma:curvesinwps2}, to deduce that up to coordinate change the
only curves of degree at most~$A^3$ which are not contained in $\{x_0 =
x_1 = 0\}$ are those listed in the third column. (In fact, for a given
family, there is usually only one of these curves up to coordinate change.)
The fourth column gives the method used to exclude the curve in question
--- usually that of Example~\ref{ex:curves_surf_1}, but in a few cases that
of Example~\ref{ex:curves_surf_2}. Each of these methods involves picking a
general surface~$T$ in some linear system with a certain base locus
containing~$C$; this linear system is given in the last column.  This
completes the proof.
\hfill $\qedsymbol$
\end{proofofcthm21s}

\vspace{1cm}
\noindent School of Mathematics \\
University of Bristol \\
Bristol $\,$ BS8 1TW\\
United Kingdom \\

\noindent Email:  \texttt{Daniel.Ryder@bristol.ac.uk}

\end{document}